\documentclass[oneside,english]{amsart}
\usepackage{lmodern}
\usepackage[T1]{fontenc}
\usepackage[latin9]{inputenc}
\usepackage{geometry}
\usepackage{mathrsfs}
\usepackage{amsbsy}
\usepackage{amstext}
\usepackage{amsthm}
\usepackage{amssymb}

\makeatletter
\numberwithin{equation}{section}
\numberwithin{figure}{section}
\theoremstyle{plain}
\newtheorem{thm}{\protect\theoremname}
\theoremstyle{remark}
\newtheorem*{rem*}{\protect\remarkname}
\theoremstyle{plain}
\newtheorem{prop}[thm]{\protect\propositionname}
\theoremstyle{plain}
\newtheorem{lem}[thm]{\protect\lemmaname}

\makeatother

\usepackage{babel}
\providecommand{\lemmaname}{Lemma}
\providecommand{\propositionname}{Proposition}
\providecommand{\remarkname}{Remark}
\providecommand{\theoremname}{Theorem}

\begin{document}
\global\long\def\e{e}%
\global\long\def\V{{\rm Vol}}%
\global\long\def\bs{\boldsymbol{\sigma}}%
\global\long\def\br{\boldsymbol{\rho}}%
\global\long\def\bp{\boldsymbol{\pi}}%
\global\long\def\btau{\boldsymbol{\tau}}%
\global\long\def\bx{\mathbf{x}}%
\global\long\def\by{\mathbf{y}}%
\global\long\def\bz{\mathbf{z}}%
\global\long\def\bv{\mathbf{v}}%
\global\long\def\bu{\mathbf{u}}%
\global\long\def\bi{\mathbf{i}}%
\global\long\def\bn{\mathbf{n}}%
\global\long\def\grad{\nabla_{sp}}%
\global\long\def\Hess{\nabla_{sp}^{2}}%
\global\long\def\lp{\Delta_{sp}}%
\global\long\def\gradE{\nabla_{\text{Euc}}}%
\global\long\def\HessE{\nabla_{\text{Euc}}^{2}}%
\global\long\def\HessEN{\hat{\nabla}_{\text{Euc}}^{2}}%
\global\long\def\ddq{\frac{d}{dR}}%
\global\long\def\qs{q_{\star}}%
\global\long\def\qss{q_{\star\star}}%
\global\long\def\lm{\lambda_{min}}%
\global\long\def\Es{E_{\star}}%
\global\long\def\As{A_{\star}}%
\global\long\def\EH{E_{\Hess}}%
\global\long\def\Esh{\hat{E}_{\star}}%
\global\long\def\ds{d_{\star}}%
\global\long\def\Cs{\mathscr{C}_{\star}}%
\global\long\def\nh{\boldsymbol{\hat{\mathbf{n}}}}%
\global\long\def\BN{\mathbb{B}^{N}}%
\global\long\def\ii{\mathbf{i}}%
\global\long\def\SN{\mathbb{S}^{N-1}}%
\global\long\def\SM{\mathbb{S}^{M-1}}%
\global\long\def\SNq{\mathbb{S}^{N-1}(q)}%
\global\long\def\SNqd{\mathbb{S}^{N-1}(q_{d})}%
\global\long\def\SNqp{\mathbb{S}^{N-1}(q_{P})}%
\global\long\def\nd{\nu^{(\delta)}}%
\global\long\def\nz{\nu^{(0)}}%
\global\long\def\cls{c_{LS}}%
\global\long\def\qls{q_{LS}}%
\global\long\def\dls{\delta_{LS}}%
\global\long\def\E{\mathbb{E}}%
\global\long\def\P{\mathbb{P}}%
\global\long\def\R{\mathbb{R}}%
\global\long\def\spp{{\rm Supp}(\mu_{P})}%
\global\long\def\indic{\mathbf{1}}%
\global\long\def\lsc{\mu_{{\rm sc}}}%
\newcommand{\SNarg}[1]{\mathbb S^{N-1}(#1)}
\global\long\def\se{s(E)}%
\global\long\def\ses{s(\Es)}%
\global\long\def\so{s(0)}%
\global\long\def\sef{s(E_{f})}%
\global\long\def\seinf{s(E_{\infty})}%
\global\long\def\L{\mathcal{L}}%
\global\long\def\gflow#1#2{\varphi_{#2}(#1)}%
\global\long\def\S{\mathscr{S}}%
\global\long\def\Frep{F^{{\rm Rep}}}%
\global\long\def\s{\mathfrak{s}}%
\global\long\def\e{e}%
\global\long\def\EsN{E_{\star,N}}%

\title{TAP approach for multi-species spherical spin glasses II: the free
energy of the pure models}
\author{Eliran Subag}
\begin{abstract}
In a companion paper we developed the generalized TAP approach for general multi-species
spherical mixed $p$-spin models. In this paper, we use it to compute the limit of the
free energy at any temperature for all pure multi-species spherical $p$-spin models, assuming that
certain free energies converge. Importantly, the pure multi-species models do not satisfy the convexity
assumption on the mixture which was crucial in the recent proofs of the Parisi formula for the
multi-species Sherrington-Kirkpatrick model by Barra et al. (2015) and Panchenko (2015) and for
the multi-species spherical mixed $p$-spin models by Bates and Sohn (2021).
\end{abstract}

\maketitle

\section{Introduction}

In the classical Sherrington-Kirkpatrick (SK) mean-field spin glass
model \cite{SK75}, the random interaction coefficients are identically
distributed for any two spins. Bipartite and, more generally, multi-species
versions of the SK model have been proposed in physics \cite{Barra,Barra14,BarraGenoveseGuerra11,KFS1,KFS2,KorenblitShender85}.
In those models the spins are divided into groups of different types,
and the strength of the interaction between any two spins depends
on their types. In this paper we will analyze the spherical $p$-spin
version of those models defined as follows. 

First, consider a finite set of species $\S$ with at least two elements,
which will be fixed throughout the paper. For each $N\geq1$, suppose
that
\[
\{1,\ldots,N\}=\bigcup_{s\in\S}I_{s},\quad\text{for some disjoint }I_{s}.
\]
The subsets $I_{s}$, of course, vary with $N$. Denoting $N_{s}=|I_{s}|$,
we will assume that the proportion of each of the species converges
\[
\lim_{N\to\infty}\frac{N_{s}}{N}=\lambda(s)\in(0,1),\quad\text{for all }s\in\S.
\]

Let $S(d)=\big\{\bx\in\R^{d}:\,\|\bx\|=\sqrt{d}\big\}$ be the sphere
of radius $\sqrt{d}$ in dimension $d$. The \emph{configuration space}
of the spherical multi-species spin glass model is
\[
S_{N}=\left\{ (\sigma_{1},\ldots,\sigma_{N})\in\R^{N}:\,\forall s\in\S,\,(\sigma_{i})_{i\in I_{s}}\in S(N_{s})\right\} .
\]
Denote $\mathbb{Z}_{+}:=\{0,1,\ldots\}$ and $|p|:=\sum_{s\in\S}p(s)$
for $p\in\mathbb{Z}_{+}^{\S}$ and define 
\[
P=\Big\{ p\in\mathbb{Z}_{+}^{\S}:\,|p|\geq1\Big\}.
\]
Given $p\in P$, define the set of indices 
\[
I(p)=\left\{ (i_{1},\ldots,i_{|p|}):\,\forall s,\,\#\{j:\,i_{j}\in I_{s}\}=p(s)\right\} ,
\]
where $\#A$ denotes the cardinality of a set $A$. 

For each $p\in P$, the spherical multi-species \emph{pure $p$-spin
Hamiltonian} is given by 
\begin{equation}
H_{N}(\bs)=H_{N,p}(\bs)=C_{N,p}\sum_{(i_{1},\ldots,i_{|p|})\in I(p)}J_{i_{1},\dots,i_{|p|}}\sigma_{i_{1}}\cdots\sigma_{i_{|p|}},\label{eq:pmodel}
\end{equation}
where $J_{i_{1},\dots,i_{k}}$ are i.i.d. standard normal variables
and 
\[
C_{N,p}^{2}=\frac{N}{\prod_{s\in\S}N_{s}^{p(s)}}\frac{\prod_{s\in\S}p(s)!}{\big(\sum_{s\in\S}p(s)\big)!}.
\]
More generally, given some nonnegative numbers $(\Delta_{p})_{p\in P}$,
the spherical multi-species \emph{mixed }$p$-spin Hamiltonian is
\begin{equation}
H_{N}(\bs)=\sum_{p\in P}\Delta_{p}H_{N,p}(\bs).\label{eq:HNmixed}
\end{equation}
Here we assume that the coefficients $J_{i_{1},\dots,i_{k}}$ are
also independent for different $k$, so that $H_{N,p}(\bs)$ are independent
for different $p$. The coefficients $\Delta_{p}$ can be encoded
in the\emph{ mixture polynomial} in the variables $x=(x(s))_{s\in\S}\in\R^{\S}$,
\[
\xi(x)=\sum_{p\in P}\Delta_{p}^{2}\prod_{s\in\S}x(s)^{p(s)}.
\]
For the mixed models which we will encounter in the current paper,
$\Delta_{p}$ is non-zero only for finitely many $p\in P$. As can
be checked by a straightforward calculation, the covariance function
of the mixed Hamiltonian $H_{N}(\bs)$ is given by 
\begin{equation}
\frac{1}{N}\E H_{N}(\bs)H_{N}(\bs')=\xi(R(\bs,\bs')),\label{eq:cov-3}
\end{equation}
where we define the overlap vector
\[
R(\bs,\bs'):=\big(R_{s}(\bs,\bs')\big)_{s\in\S},\quad R_{s}(\bs,\bs'):=N_{s}^{-1}\sum_{i\in I_{s}}\sigma_{i}\sigma_{i}'.
\]

Identifying $S_{N}$ with the product space $\prod_{s\in\S}S(N_{s})$,
let $\mu$ be the product of the uniform measures on each of the spheres
$S(N_{s})$. The (mean) free energy $F_{N}(\beta)$ and partition
function $Z_{N}(\beta)$ at inverse-temperature $\beta>0$ are defined
by 
\begin{equation}
F_{N}(\beta):=\frac{1}{N}\E\log Z_{N,\beta}:=\frac{1}{N}\E\log\int_{S_{N}}e^{\beta H_{N}(\bs)}d\mu(\bs).\label{eq:F-2}
\end{equation}
The (mean) normalized ground state energy is
\begin{equation}
\EsN:=\frac{1}{N}\E\max_{\bs\in S_{N}}H_{N}(\bs).\label{eq:Esq-1-1}
\end{equation}
The Gibbs measure at inverse-temperature $\beta>0$ is the random
probability measure on $S_{N}$ with density
\[
\frac{dG_{N,\beta}}{d\mu}(\bs)=Z_{N,\beta}^{-1}e^{\beta H_{N}(\bs)}.
\]

Very recently, Bates and Sohn \cite{BatesSohn1,BatesSohn2} proved
a Parisi formula \cite{ParisiFormula,Parisi1980} for spherical multi-species
mixed $p$-spin models such that the mixture $\xi(x)$ is a convex
function. (For the proof of the formula in the single-species case,
see \cite{Chen,Panch,Talag,Talag2}.) The same convexity condition
was required in an earlier work of Barra, Contucci, Mingione and Tantari
\cite{Barra}, where they proved for the multi-species Sherrington-Kirkpatrick
model that a Parisi type formula is an upper bound for the limiting
free energy using an analogue of Guerra's interpolation \cite{GuerraBound}.
Panchenko proved the matching lower bound in \cite{PanchenkoMulti}
using an analogue of the Aizenman-Sims-Starr scheme \cite{ASSscheme}.
One of the main novelties in his paper was the proof that the overlaps
$R_{s}(\bs,\bs')$ of independent samples from the Gibbs measure are
synchronized for different species. The synchronization mechanism
of \cite{PanchenkoMulti} was also used in the spherical setting in
\cite{BatesSohn1}. The combination of \cite{Barra} and \cite{PanchenkoMulti}
establishes the Parisi formula for the multi-species Sherrington-Kirkpatrick
model, assuming the convexity of $\xi(x)$.

In a companion paper \cite{TAPmulti1} we developed the generalized
TAP approach for general multi-species spherical mixed $p$-spin models.
In this paper we use it to compute the free energy of the \emph{pure}
$p$-spin spherical multi-species models with Hamiltonian (\ref{eq:pmodel}).
For the pure $p$-spin multi-species mixture $\xi(x)=\prod_{s\in\S}x(s)^{p(s)}$,
the Hessian is given by 
\[
\nabla^{2}\xi(x):=\Big(\frac{d}{dx(s)}\frac{d}{dx(t)}\xi(x)\Big)_{s,t\in\S}=\xi(x)\bigg[\Big(\frac{p(s)}{x(s)}\Big)_{s\in\S}\cdot\Big(\frac{p(s)}{x(s)}\Big)_{s\in\S}^{T}-\text{diag}\Big(\frac{p(s)}{x(s)^{2}}\Big)_{s\in\S}\bigg],
\]
assuming that $x(s)\neq0$ for all $s\in\S$ and WLOG that $\S=\{1,\ldots,|\S|\}$,
so that the matrix operations above make sense. Thus, for any vector
$v$ orthogonal to $\big(\frac{p(s)}{x(s)}\big)_{s\in\S}$, assuming
that $p(s)\geq1$ and $x(s)>0$,
\[
v^{T}\nabla^{2}\xi(x)v<0.
\]
From this, it is easy to verify that $\nabla^{2}\xi(x)$ has one positive
eigenvalue and $|\S|-1$ negative eigenvalues and $\xi(x)$ is therefore
not convex anywhere in $[0,1]^{\S}$.

In Section \ref{sec:TAP}, we will define certain mixtures $\xi_{q}(x)$
(see (\ref{eq:xiq}) below) for any $q\in[0,1)^{\S}$. To those mixtures
we associate the free energies $F_{N}^{q}(\beta)$, (see (\ref{eq:FqN}))
defined similarly to (\ref{eq:F-2}). In our main results we will
assume that they converge as $N\to\infty$. Namely, for some $F(\beta,q)$,
\begin{equation}
\forall q\in[0,1)^{\S},\,\beta>0:\quad\lim_{N\to\infty}F_{N}^{q}(\beta)=F(\beta,q).\label{eq:convergence}
\end{equation}
This is currently only proven \cite{BatesSohn1} if $\xi_{q}(x)$
is convex on $[0,1]^{\S}$, but is obviously expected to be true in
general. Since for $q=0$ the mixture coincides with that of the original
model $\xi_{q}(x)=\xi(x)$, under this assumption the limit of the
free energy of the pure model exists
\begin{equation}
\lim_{N\to\infty}F_{N}(\beta)=F(\beta):=F(\beta,0).\label{eq:Fconv}
\end{equation}
One can also easily check\footnote{Using bounds on the Lipschitz constant of $H_{N}(\bs)$ as in \cite[Lemma 25]{TAPmulti1}
 and the fact that $F(\beta)$ is convex.} that the ground-state energy also has a limit
\[
\Es:=\lim_{N\to\infty}\EsN=\lim_{\beta\to\infty}\frac{1}{\beta}F(\beta).
\]

From Jensen's inequality,
\[
F(\beta)\leq\frac{1}{2}\beta^{2}\xi(1),
\]
where for $c\in\R$ we write $\xi(c)$ for the evaluation of $\xi$
at the constant vector $x\equiv c$. It is not difficult to check
that if $|p|\geq2$, for very small $\beta\geq0$,
\[
\lim_{N\to\infty}\frac{1}{N}\log\frac{\E\big(Z_{N,\beta}^{2}\big)}{\big(\E Z_{N,\beta}\big)^{2}}=0.
\]
Therefore, in this case, by the Paley\textendash Zygmund inequality
and concentration of the free energy, $F(\beta)=\frac{1}{2}\beta^{2}\xi(1)$.
Also, if the latter equality holds for some $\beta,$ then it also
holds for any $0<\beta'<\beta$.\footnote{\label{fn:sub_crt}In distribution, $\beta H_{N}(\bs)=\beta'H_{N}^{1}(\bs)+\sqrt{\beta^{2}-\beta'^{2}}H_{N}^{2}(\bs)$
for independent copies $H_{N}^{i}(\bs)$ of $H_{N}(\bs)$. By conditioning
on $H_{N}^{1}(\bs)$ and using Jensen's inequality, $F_{N,\beta}\leq F_{N,\beta'}+\frac{1}{2}(\beta^{2}-\beta'^{2})\xi(1)$.} Hence, there exists some critical $\beta_{c}>0$ such that
\begin{equation}
F(\beta)=\frac{1}{2}\beta^{2}\xi(1)\iff\beta\leq\beta_{c}.\label{eq:bcdef}
\end{equation}

Let $G_{N,\beta}^{\otimes n}$ be the $n$-fold product measure of
the Gibbs measure with itself. In \cite{FElandscape,TAPmulti1} we
introduced the following notion of multi-samplable overlaps. We say
that $q\in[0,1)^{\S}$ is multi-samplable at inverse-temperature $\beta$
if for any $n\geq1$ and $\epsilon>0$, 
\begin{equation}
\lim_{N\to\infty}\frac{1}{N}\log\E G_{N,\beta}^{\otimes n}\left\{ \forall i<j\leq n,\,s\in\S:\,\big|R_{s}(\bs^{i},\bs^{j})-q(s)\big|<\epsilon\right\} =0.\label{eq:multisamp-1-1}
\end{equation}
We will say that a multi-samplable overlap $q$ is maximal if there
is no other multi-samplable overlap $q'\neq q$ such that $q'(s)\geq q(s)$
for all $s\in\S$. 

Without loss of generality, we will henceforth assume that $p(s)\geq1$
for any $s\in\S$, since otherwise we can integrate out the corresponding
spins. Let
\begin{align*}
V(q) & :=-\sum_{t\in\S}\lambda(t)\log(1-q(t)),\\
U(q) & :=1+\sum_{t\in\S}(1-q(t))\frac{p(t)}{q(t)},
\end{align*}
and define
\[
\Phi(q):=\frac{V(q)}{U(q)}\qquad\text{and}\qquad\Omega(q):=V(q)U(q).
\]
The following are our main results. 
\begin{thm}
[Ground state energy and critical parameters]\label{thm:main1}For
any pure $p$-spin model $\xi(x)=\prod_{s\in\S}x(s)^{p(s)}$ with
$|p|\geq3$, assuming the convergence (\ref{eq:convergence}), the
following hold:
\begin{enumerate}
\item \label{enu:thm1pt1}At the critical inverse-temperature $\beta_{c}$,
there is exactly one maximal multi-samplable overlap $q_{c}$ and
it is equal to the unique solution $q\in(0,1)^{\S}$ of the system
of equations
\begin{equation}
\forall s\in\S:\quad\frac{\lambda(s)}{p(s)}\frac{q(s)^{2}}{1-q(s)}=\Phi(q).\label{eq:qc_formula}
\end{equation}
\item \label{enu:thm1pt2}The critical inverse-temperature is given by
\begin{equation}
\beta_{c}=\sqrt{\frac{\Phi(q_{c})}{\xi(q_{c})}}.\label{eq:bc_formula}
\end{equation}
\item \label{enu:thm1pt3}The ground state energy is given by 
\begin{equation}
\Es=\sqrt{\Omega(q_{c})}.\label{eq:Es_formula}
\end{equation}
\end{enumerate}
\end{thm}

\begin{rem*}
In Section \ref{subsec:Proof-of-Lemma_uniqqc} we will see that for
some function $z\mapsto q_{z}$ from $(0,\infty)$ to $(0,1)^{\S}$,
the fixed point equation $z=\Phi(q_{z})$ has a unique solution in
$(0,\infty)$ and for this solution $z$, the vector $q=q_{z}$ is
the unique solution of the system of equations (\ref{eq:qc_formula})
in $|\S|$ variables.
\end{rem*}

\begin{thm}
[Free energy]\label{thm:main2} For any pure $p$-spin model $\xi(x)=\prod_{s\in\S}x(s)^{p(s)}$
with $|p|\geq3$, assuming the convergence (\ref{eq:convergence}),
for any $\beta>\beta_{c}$, there is exactly one maximal multi-samplable
$q$ and it is defined by the system of equations
\begin{equation}
\forall s\in\S:\quad\frac{1-q(s)}{q(s)}=\frac{-\sqrt{\vphantom{\frac{\lambda(s)}{p(s)}}\Phi(q_{c})}+\sqrt{\Phi(q_{c})+4\frac{\lambda(s)}{p(s)}}}{2y_{\star}(\beta)},\label{eq:q}
\end{equation}
where $y_{\star}(\beta)$ is the larger of the two solutions in $(0,\infty)$
of 
\begin{equation}
y^{2}\prod_{s\in\S}\left(\frac{-\sqrt{\vphantom{\frac{\lambda(s)}{p(s)}}\Phi(q_{c})}+\sqrt{\Phi(q_{c})+4\frac{\lambda(s)}{p(s)}}}{2y}+1\right)^{p(s)}=\beta^{2}\label{eq:qb}
\end{equation}
and $q_{c}$ is the maximal multi-samplable at $\beta_{c}$ from Theorem
\ref{thm:main1}. With the same overlap $q$, the free energy is given
by 
\begin{equation}
\begin{aligned}F(\beta) & =\beta\sqrt{\xi(q)}\Es+\frac{1}{2}\sum_{s\in\S}\lambda(s)\log(1-q(s))+\frac{1}{2}\beta^{2}\left(\xi(1)-\xi(q)-\xi(q)\sum_{s\in\S}p(s)\frac{1-q(s)}{q(s)}\right).\end{aligned}
\label{eq:FThm}
\end{equation}
\end{thm}

The only multi-species pure $p$-spin model that is not covered by
the above is the case where $|p|=2$, and thus, $|\S|=2$ with $p(s)=1$
for both species. This is the bipartite model studied by Baik and
Lee in \cite{BaikLeeBipartite}. We analyze this model in the Appendix,
where we compute the free energy and maximal multi-samplable overlap
for all $\beta$ from the TAP representation. 

The theorems above generalize to the multi-species case the main results
for the single-species pure models in \cite{SubagTAPpspin}, derived
using the TAP representation from \cite{FElandscape}. For the mixed
single-species spherical models, the free energy is given by the well-known
Crisanti-Sommers representation \cite{Crisanti1992} of the Parisi
formula \cite{Parisi,ParisiFormula}, first prove by Talagrand \cite{Talag}
and extended by Chen \cite{Chen}. Several earlier results are connected
to the TAP representation in the single-species case. Belius and Kistler
\cite{BeliusKistlerTAP} established a TAP representation for the
free energy for the spherical $2$-spin model with an external field.
In \cite{geometryGibbs} the free energy was computed for single-species
spherical pure $p$-spin models at low enough temperature, building
on results about critical points \cite{A-BA-C,ABA2,2nd,pspinext},
and was shown to be given by the TAP representation. A similar result
was proved for mixed models close to pure by Ben Arous, Zeitouni and
the author \cite{geometryMixed}. For more recent results on the TAP
equations and representation for the single-species models with Ising
spins see \cite{AuffingerJagannthSpinDist,AuffingerJagannathTAP,BolthausenTAP,BolthausenMorita,ChatterjeeTAP,ChenPanchenkoTAP,TalagrandBookI}.

Multi-species models with two species, i.e., $\S=\{s_{1},s_{2}\}$,
are called bipartite models. Baik and Lee \cite{BaikLeeBipartite}
computed the free energy and the limiting law of its fluctuations
for the pure $p$-spin spherical bipartite model with $p(s_{1})=p(s_{2})=1$
using tools from random matrix theory. Auffinger and Chen \cite{AuffingerChenBipartite}
proved that for mixed spherical bipartite models and in the presence
of an external field, if $\xi(1)$ and the strength of the field are
small enough, then the limiting free energy is given by the replica
symmetric solution of an analogue of the Crisanti-Sommers formula
\cite{Crisanti1992}. Certain bounds on the exponential growth rate
of  the mean number of minimum points, known as `complexity', were
also derived in \cite{AuffingerChenBipartite}. McKenna \cite{McKennaComplexity}
computed the exact asymptotics for the mean complexity of minimum
and critical points. For the pure bipartite model, Kivimae \cite{Kivimae}
proved that for large energies, the second moment of the complexity
matches its mean squared at exponential scale, assuming $p(s_{1})$
and $p(s_{2})$ are large enough. Combining \cite{McKennaComplexity}
and \cite{Kivimae} yields a variational formula for the ground state
energy of the pure bipartite models. As mentioned above, Bates and
Sohn proved a Parisi formula for the multi-species spherical mixed
$p$-spin models in \cite{BatesSohn1}, assuming the convexity of
$\xi(x)$. In \cite{BatesSohn2} they showed it admits a representation
analogous to the Crisanti-Sommers formula \cite{Crisanti1992}. 

We also mentioned above the results of Barra, Contucci, Mingione and
Tantari \cite{Barra} and Panchenko \cite{PanchenkoMulti} which combined
prove the Parisi formula for the multi-species SK model, with the
same convexity assumption on $\xi(x)$. Bates, Sloman and Sohn \cite{BatesSlomanSohn}
study the replica symmetric phase of the latter formula. Mourrat derived
in \cite{MourratBipartite} an upper bound in terms of an infinite-dimensional Hamilton-Jacobi equation for the free energy of
the bipartite SK model with only inter-species interactions, whose
mixture is non-convex. In \cite{MourratVector}, he generalized the bound for essentially arbitrary fully-connected models.

In the next section we will recall the generalized TAP representation
from \cite{TAPmulti1}, and briefly explain how it will be used to
obtain a system of equations involving the various quantities from
Theorems \ref{thm:main1} and \ref{thm:main2}. In Section \ref{sec:Additional-equations},
we will prove the latter equations, stated in Proposition \ref{prop:eqns}
below. In Section \ref{sec:PfThm1} we will prove Theorem \ref{thm:main1}
concerning the critical inverse-temperature and ground state energy.
In Section \ref{sec:pf_Thm2} we will prove Theorem \ref{thm:main2}
concerning the value of free energy for super-critical inverse-temperatures
$\beta$.

\section{\label{sec:TAP}The generalized TAP approach}

In the 70s, Thouless, Anderson and Palmer \cite{TAP} invented their
celebrated approach to analyze the SK model. Their approach was further
developed in physics, see e.g. \cite{TAP-SK1,TAP-SK3,TAP-pSPSG1,TAP-SK2,TAP-pS-Ising3,KurchanParisiVirasoro,Plefka},
with the general idea that for large $N$, $F_{N}(\beta)$ is approximated
by the free energies associated to the `physical' solutions of the
TAP equations. In particular, the single-species spherical pure $p$-spin
models were analyzed non-rigorously using the TAP approach by Kurchan,
Parisi and Virasoro in \cite{KurchanParisiVirasoro} and Crisanti
and Sommers in \cite{TAP-pSPSG1}.

Recently, we developed in \cite{FElandscape} a generalized TAP approach
for the single-species spherical models, which can be applied to any
multi-samplable overlap. The approach was extended to mixed models
with Ising spins by Chen, Panchenko and the author \cite{TAPChenPanchenkoSubag,TAPIIChenPanchenkoSubag}.
In a companion paper \cite{TAPmulti1} we developed the generalized
TAP approach for general multi-species spherical mixed $p$-spin glass
models. The proofs of our main theorems will be based on the TAP representation
for the free energy from \cite{TAPmulti1}. In this section we recall
the representation in the setting of the pure models, and explain
how it will be used to obtain a system of equations involving the
various quantities from Theorems \ref{thm:main1} and \ref{thm:main2}.

Consider the pure $p$-spin mixture $\xi(x)=\prod_{s\in\S}x(s)^{p(s)}$
and let $q\in[0,1)^{\S}$. Define the mixture 
\begin{equation}
\xi_{q}(x)=\xi((1-q)x+q)-\xi(q)-\sum_{s\in\S}\big((1-q)\nabla\xi(q)x\big)(s),\label{eq:xiq}
\end{equation}
where 
\[
\nabla\xi(q)(s):=\frac{d}{dq(s)}\xi(q)=\frac{p(s)}{q(s)}\xi(q).
\]
Here and in the sequel, for $q(s)=0$, we use the convention that division
by $q(s)$ should be interpreted by cancellation with the corresponding
power of $q(s)$ in $\xi(q)$, e.g., $\xi(q)/q(s)=q(s)^{p(s)-1}\prod_{t\neq s}q(t)^{p(t)}$.
Also, in (\ref{eq:xiq}) for any two functions $f$ and $g$ of $s\in\S$,
we use the usual notation for the sum $(f+g)(s)=f(s)+g(s)$ and product
$(fg)(s)=f(s)g(s)$. 

One can check that 
\[
\xi_{q}(x)=\sum_{k\in P:\,|k|\geq2,\,k\leq p}\Delta_{q,k}^{2}\prod_{s\in\S}x(s)^{k(s)},
\]
where we write $k\leq p$ if $k(s)\leq p(s)$ for all $s\in\S$ and
\begin{equation}
\Delta_{q,k}^{2}:=\prod_{s\in\S}\binom{p(s)}{k(s)}(1-q(s))^{k(s)}q(s)^{p(s)-k(s)}.\label{eq:Deltaqp}
\end{equation}
Denote the Hamiltonian corresponding to the mixture $\xi_{q}(x)$
by $H_{N}^{q}(\bs)$ and its free energy by
\begin{equation}
F_{N}^{q}(\beta):=\frac{1}{N}\E\log\int_{S_{N}}e^{\beta H_{N}^{q}(\bs)}d\mu(\bs).\label{eq:FqN}
\end{equation}
Recall that we assume in (\ref{eq:convergence}) that the latter free
energies converge to some $F(\beta,q)$.

The mixture $\xi_{q}(x)$ has the following meaning. For $\bs$ in
\begin{equation}
\Big\{\bs\in S_{N}:\,\forall s\in\S,\,R_{s}(\bs,m)=R_{s}(m,m)\Big\}\label{eq:B0}
\end{equation}
define $\tilde{\bs}$ by rescaling $\bs-m$ back to $S_{N}$,
\begin{equation}
\tilde{\sigma}_{i}:=\sqrt{\frac{1}{1-q(s)}}(\sigma_{i}-m_{i}),\quad\text{if }i\in I_{s}.\label{eq:sigmatilde}
\end{equation}
On the image of (\ref{eq:B0}) by the mapping $\bs\mapsto\tilde{\bs}$,
define Hamiltonian $\tilde{H}_{N}(\tilde{\bs})=H_{N}(\bs)-H_{N}(m)$.
Then, for any two points $\bs^{1},\,\bs^{2}$ in (\ref{eq:B0}),
\[
\frac{1}{N}\E\tilde{H}_{N}(\tilde{\bs}^{1})\tilde{H}_{N}(\tilde{\bs}^{2})=\xi_{q}(R(\tilde{\bs}^{1},\tilde{\bs}^{2}))+\sum_{s\in\S}\Delta_{q,s}^{2}R_{s}(\tilde{\bs}^{1},\tilde{\bs}^{2}),
\]
where by abuse of notation we use $\Delta_{q,s}$ to denote the coefficient
as in (\ref{eq:Deltaqp}) for $k$ such that $k(s)=1$ and $k(t)=0$
for any $t\neq s$.

Specialized to the pure $p$-spin model, the TAP representation for
the free energy proved in \cite{TAPmulti1} states that if $q\in[0,1)^{\S}$
is multi-samplable then
\begin{equation}
F(\beta)=\beta\sqrt{\xi(q)}\Es+\frac{1}{2}\sum_{s\in\S}\lambda(s)\log(1-q(s))+F(\beta,q),\label{eq:TAPrep-1}
\end{equation}
and if $q\in[0,1)^{\S}$ is not multi-samplable 
\begin{equation}
F(\beta)>\beta\sqrt{\xi(q)}\Es+\frac{1}{2}\sum_{s\in\S}\lambda(s)\log(1-q(s))+F(\beta,q).\label{eq:TAPineq}
\end{equation}
Here we used the fact that
\begin{equation}
\EsN(q):=\frac{1}{N}\E\max_{m\in S_{N}(q)}H_{N}(m)\longrightarrow\sqrt{\xi(q)}\Es,\label{eq:Eqconv}
\end{equation}
where 
\[
S_{N}(q):=\left\{ m\in\R^{N}:\,R(m,m)=q\right\} .
\]
Moreover, it was proved in \cite{TAPmulti1} that if $q\in[0,1)^{\S}$
is a maximal multi-samplable overlap, (assuming the convergence in
(\ref{eq:convergence})) then
\begin{equation}
F(\beta,q)=\frac{1}{2}\beta^{2}\xi_{q}(1).\label{eq:Onsager}
\end{equation}

Suppose that $q_{c}$ is a maximal multi-samplable overlap at the
critical inverse-temperature $\beta_{c}$. For the critical inverse-temperature,
the free energy is equal to $F(\beta_{c})=\frac{1}{2}\beta_{c}^{2}\xi(1)$.
Combining this with (\ref{eq:TAPrep-1}) and (\ref{eq:Onsager}) gives
us one equation for the triplet $(\beta_{c},q_{c},\Es)$. 

To obtain additional equations, we will exploit the fact that for
overlaps which are not multi-samplable we have the inequality (\ref{eq:TAPineq}).
Since $\beta_{c}>0$, this implies that the derivatives in $\beta$
of both sides of the TAP representation (\ref{eq:TAPrep-1}) at $\beta=\beta_{c}$
and $q=q_{c}$ must be equal, provided they exist. It is simple to
see that $F'(\beta_{c})=\beta_{c}\xi(1)$. We will also see that $\frac{d}{d\beta}F(\beta,q)=\beta\xi_{q}(1)$,
if $q$ is a maximal multi-samplable overlap at $\beta$ (see Lemma
\ref{lem:derivative_bc}). This will give us an additional equation
in $(\beta_{c},q_{c},\Es)$. 

Finally, we will show that $q_{c}(s)>0$, and thus by a similar argument,
the derivatives in $q(s)$ of both sides of (\ref{eq:TAPrep-1}) are
equal. The left-hand side does not depend on $q$ and its derivative
is therefore $0$. We will show that $\frac{d}{dq(s)}F(\beta,q)=\frac{1}{2}\beta^{2}\frac{d}{dq(s)}\xi_{q}(1)$,
if $q$ is a maximal multi-samplable overlap at $\beta$ (see Lemma
\ref{lem:derivative_qb}). This will give us another equation for
each $s\in\S$. In total, the number of equations we will have is
exactly the number of scalar variables $|\S|+2$ (since $q_{c}$ is
a vector of length $|\S|$). Precisely, we will prove the following
proposition in the next section. 
\begin{prop}
\label{prop:eqns}For the pure $p$-spin model $\xi(q)=\prod_{s\in\S}q(s)^{p(s)}$
with $|p|\geq3$, assuming (\ref{eq:convergence}), if $q_{c}$ is
a maximal multi-samplable overlap at the critical inverse-temperature
$\beta=\beta_{c}$, then for all $s\in\S$, $q_{c}(s)>0$ and the
triple $(\beta_{c},q_{c},\Es)$ solves 
\begin{align}
\frac{\lambda(s)}{1-q(s)}+\beta^{2}\xi(q)\left(-p(s)\frac{1-q(s)}{q(s)^{2}}+\frac{p(s)}{q(s)}\sum_{t\in\S}p(t)\frac{1-q(t)}{q(t)}\right) & =\beta E\sqrt{\xi(q)}\frac{p(s)}{q(s)},\label{eq:I}\\
\beta^{2}\xi(q)\left(1+\sum_{t\in\S}p(t)\frac{1-q(t)}{q(t)}\right)=-\sum_{t\in\S}\lambda(t)\log(1-q(t)) & =\beta E\sqrt{\xi(q)}.\label{eq:II}
\end{align}
Moreover, for any $\beta>\beta_{c}$, if $q_{\beta}$ is a maximal
multi-samplable overlap then, for all $s\in\S$, $q_{\beta}(s)>0$
and $(\beta,q_{\beta},\Es)$ solves (\ref{eq:I}).
\end{prop}

For $\beta>\beta_{c}$, the proposition gives us the $|\S|$ equations
(\ref{eq:I}) derived from equality of derivatives in $q(s)$. The
two other equalities (\ref{eq:II}), rely on the fact that at $\beta=\beta_{c}$,
$F(\beta_{c})=\frac{1}{2}\beta_{c}^{2}\xi(1)$ and $F'(\beta_{c})=\beta_{c}\xi(1)$
and are therefore not valid for $\beta>\beta_{c}$. However, for such
$\beta$ we will only need to compute the maximal multi-samplable
overlap $q$, and we will already have the value of $\Es$ from Theorem
\ref{thm:main1}. Thus, in the setting of Theorem \ref{thm:main2},
we will need to solve the $|\S|$ equations (\ref{eq:I}) for the
same number of variables $q(s)$. Once we do that, we will be able
to express the free energy by (\ref{eq:TAPrep-1}) to obtain (\ref{eq:FThm}).

\section{\label{sec:Additional-equations}Additional equations from the TAP
representation: Proof of Proposition \ref{prop:eqns}}

In this section we prove Proposition \ref{prop:eqns} from the equality
of derivatives as explained in the previous section. We will prove
the following two lemmas in Sections \ref{subsec:pfddq} and \ref{subsec:pfddb}.
Below, the existence of the derivatives is part of the statement.
As mentioned above, to lighten the notation we use the convention
that if $q(t)=0$, one should interpret the division by $q(t)$ by
cancellation with corresponding the power of $q(t)$ in $\xi(q)$,
e.g., $\xi(q)/q(t)=q(t)^{p(t)-1}\prod_{s\neq t}q(s)^{p(s)}$.
\begin{lem}
\label{lem:derivative_qb}For any pure $p$-spin model and any $\beta\geq0$,
if $q$ is a maximal multi-samplable overlap, then for any $s\in\S$
such that $q(s)>0$, $(\beta,q,\Es)$ solves the equations 
\begin{align}
\frac{d}{dq(s)}F(\beta,q) & =\frac{d}{dq(s)}\frac{1}{2}\beta^{2}\xi_{q}(1)\label{eq:ddqF}\\
 & =\frac{1}{2}\beta^{2}\xi(q)\left(p(s)\frac{1-q(s)}{q(s)^{2}}-\frac{p(s)}{q(s)}\sum_{t\in\S}p(t)\frac{1-q(t)}{q(t)}\right),\nonumber \\
0 & =\frac{1}{2}\beta E\sqrt{\xi(q)}\frac{p(s)}{q(s)}-\frac{1}{2}\frac{\lambda(s)}{1-q(s)}\label{eq:ddqF2}\\
 & +\frac{1}{2}\beta^{2}\xi(q)\left(p(s)\frac{1-q(s)}{q(s)^{2}}-\frac{p(s)}{q(s)}\sum_{t\in\S}p(t)\frac{1-q(t)}{q(t)}\right).\nonumber 
\end{align}
\end{lem}

\begin{lem}
\label{lem:derivative_bc}For any pure $p$-spin model and $\beta>0$,
if $q$ is a maximal multi-samplable overlap then 
\begin{equation}
\frac{d}{d\beta}F(\beta,q)=\frac{d}{d\beta}\frac{1}{2}\beta^{2}\xi_{q}(1)=\beta\left(\xi(1)-\xi(q)-\xi(q)\sum_{t\in\S}p(t)\frac{1-q(t)}{q(t)}\right).\label{eq:ddb-2}
\end{equation}
For $\beta=\beta_{c}$, we further have that
\begin{equation}
\beta\xi(1)=\Es\sqrt{\xi(q)}+\beta\left(\xi(1)-\xi(q)-\xi(q)\sum_{t\in\S}p(t)\frac{1-q(t)}{q(t)}\right).\label{eq:ddb2}
\end{equation}
\end{lem}

We will prove the following lemma in Section \ref{subsec:q_nonzero}.
Note that here the assumption that $p(s)\geq1$ for all $s\in\S$
is essential.
\begin{lem}
\label{lem:q_nonzero}For any pure $p$-spin model, if $q\in[0,1)^{\S}$
is a maximal multi-samplable, then either $q(s)=0$ for all $s\in\S$
or $q(s)>0$ for all $s\in\S$.
\end{lem}

Suppose that $q$ is a maximal multi-samplable overlap at inverse-temperature
$\beta$. From the TAP representation (\ref{eq:TAPrep-1}) and (\ref{eq:Onsager}),
we have that if $q(s)=0$ for all $s\in\S$, then
\[
F(\beta)=\frac{1}{2}\beta^{2}\xi_{q}(1)=\frac{1}{2}\beta^{2}\xi(1),
\]
that is, $\beta\leq\beta_{c}$. Hence, by the lemma above, for any
$\beta>\beta_{c}$, if $q$ is a maximal multi-samplable overlap then
$q(s)>0$ for all $s\in\S$. In Section \ref{subsec:pfq>0}, we will
prove the following lemma by which the same is true at the critical
$\beta=\beta_{c}$.
\begin{lem}
\label{lem:q>0}For any pure $p$-spin model with $|p|\geq3$, if
$q$ is a maximal multi-samplable overlap at $\beta_{c}$, then $q(s)>0$
for all $s\in\S$.
\end{lem}

From the above, for any $\beta\geq\beta_{c}$, if $q$ is a maximal
multi-samplable overlap then $q(s)>0$ for all $s\in\S$. Hence, such
$\beta$ and $q$ satisfy the equations of Lemmas \ref{lem:derivative_qb}
and \ref{lem:derivative_bc}. Moreover, from the TAP representation
(\ref{eq:TAPrep-1}) and (\ref{eq:Onsager}),
\begin{equation}
\begin{aligned}F(\beta) & =\beta\sqrt{\xi(q)}\Es+\frac{1}{2}\sum_{s\in\S}\lambda(s)\log(1-q(s))\\
 & +\frac{1}{2}\beta^{2}\left(\xi(1)-\xi(q)-\xi(q)\sum_{s\in\S}p(s)\frac{1-q(s)}{q(s)}\right).
\end{aligned}
\label{eq:TAPq>0}
\end{equation}

Equation (\ref{eq:I}) follows from (\ref{eq:ddqF2}). For $\beta=\beta_{c}$,
the equality of the right-hand and left-hand sides of (\ref{eq:II})
follows from (\ref{eq:ddb2}). Lastly, for $\beta=\beta_{c}$, from
(\ref{eq:TAPq>0}) and the fact $F(\beta)=\frac{1}{2}\beta^{2}\xi(1)$
we have that 
\[
\beta\sqrt{\xi(q)}\Es+\frac{1}{2}\sum_{s\in\S}\lambda(s)\log(1-q(s))=\frac{1}{2}\beta^{2}\xi(q)\left(1+\sum_{s\in\S}p(s)\frac{1-q(s)}{q(s)}\right).
\]
Combined with the equality of the right-hand and left-hand sides of
(\ref{eq:II}), this completes the proof of (\ref{eq:II}).\qed 

It remains to prove the lemmas above.

\subsection{\label{subsec:pfddq}Proof of Lemma \ref{lem:derivative_qb}}

Let $\beta\geq0$ and suppose that $q$ is a maximal multi-samplable
overlap. Assume that $q(s)>0$ for some $s\in\S$. By Corollary 18
of \cite{TAPmulti1}, if we denote by $G_{N,\beta,q}$ the Gibbs measure
corresponding to the mixture $\xi_{q}(x)$, then, for any $\epsilon>0$,
\begin{equation}
\lim_{N\to\infty}\E G_{N,\beta,q}^{\otimes2}\big\{\exists s\in\S:\,|R_{s}(\bs^{1},\bs^{2})|>\epsilon\big\}=0.\label{eq:RS_Gibbs}
\end{equation}
Using this, by the same proof as that of (2.5) in Lemma 5 of \cite{SubagTAPpspin},
the derivative $\frac{d}{dq(s)}F(\beta,q)$ exists and satisfies (\ref{eq:ddqF}).

As explained in the beginning of the section, since the TAP representation
holds as an inequality for any $q\in[0,1)^{\S}$ and as equality for
multi-samplable $q$, the derivative in $q(s)$ of the right-hand
side of (\ref{eq:TAPrep-1}) is equal to that of the left-hand side,
namely,  to zero. Hence, (\ref{eq:ddqF2}) follows from (\ref{eq:ddqF}).\qed

\subsection{\label{subsec:pfddb}Proof of Lemma \ref{lem:derivative_bc}}

Recall that for any $\beta\leq\beta_{c}$, $F(\beta)=\frac{1}{2}\beta^{2}\xi(1)$.
Hence, the derivative of the free energy from the left at $\beta_{c}$
is equal to $\frac{d}{d\beta}^{-}F(\beta_{c})=\beta_{c}\xi(1)$. By
H\"{o}lder's inequality $F(\beta)$ is convex, and therefore the
derivative from the right exists and satisfies $\frac{d}{d\beta}^{-}F(\beta_{c})\leq\frac{d}{d\beta}^{+}F(\beta_{c})$.
By Jensen's inequality, for any $\beta>\beta_{c}$, $F(\beta)\leq\frac{1}{2}\beta^{2}\xi(1)$.
Hence, $\frac{d}{d\beta}^{+}F(\beta_{c})\leq\beta_{c}\xi(1)$. That
is, both one-sided derivatives exist and are equal and thus $\frac{d}{d\beta}F(\beta_{c})=\beta_{c}\xi(1)$.

Let $\beta$ and suppose that $q$ is a maximal multi-samplable overlap.
Then, (\ref{eq:Onsager}) exactly means that $\beta$ is less than
or equal to the critical inverse-temperature of the model with mixture
$\xi_{q}(x)$. If $\beta$ is equal to the critical inverse-temperature
of $\xi_{q}(x)$, then applying the argument above to the corresponding
free energy $F(\beta,q)$ implies (\ref{eq:ddb-2}). If $\beta$ is
strictly less than the critical inverse-temperature of $\xi_{q}(x)$,
then for $\beta_{0}$ in a neighborhood of $\beta$, $F(\beta_{0},q)=\frac{1}{2}\beta_{0}^{2}\xi_{q}(1)$
and (\ref{eq:ddb-2}) follows.

With the maximal multi-samplable $q_{c}$, the TAP representation
(\ref{eq:TAPrep-1}) holds as an inequality $\geq$ for any $\beta$
and as equality for $\beta=\beta_{c}$. Hence, the derivative in $\beta$
of both sides of (\ref{eq:TAPrep-1}) is equal at $\beta_{c}$. From
this, (\ref{eq:ddb2}) follows from (\ref{eq:ddb-2}) and the fact
that $\frac{d}{d\beta}F(\beta_{c})=\beta_{c}\xi(1)$.\qed

\subsection{\label{subsec:q_nonzero}Proof of Lemma \ref{lem:q_nonzero}}

Suppose that $q$ is a maximal multi-samplable overlap and assume
towards contradiction that for some $s,\,s'\in\S$, $q(s)>0$ and
$q(s')=0$. By Lemma \ref{lem:derivative_qb},
\begin{equation}
\beta\Es\sqrt{\xi(q)}\frac{p(s)}{q(s)}-\frac{\lambda(s)}{1-q(s)}+\beta^{2}\xi(q)\left(p(s)\frac{1-q(s)}{q(s)^{2}}-\frac{p(s)}{q(s)}\sum_{t\in\S}p(t)\frac{1-q(t)}{q(t)}\right)=0,\label{eq:20051}
\end{equation}
where we use the same notation convention mentioned before the same
lemma. Note that since $q(s')=0$, $\xi(q)=0$ and the first term
above is zero. The second term is strictly negative, and the third
in non-positive. Hence, we arrive at a contradiction, and the proof
is completed.\qed

\subsection{\label{subsec:pfq>0}Proof of Lemma \ref{lem:q>0} }

Recall that $\beta_{c}>0$ and let $\beta_{k}\in(\beta_{c},2\beta_{c})$
be a decreasing sequence that converges to $\beta_{c}$. For any $k$,
let $q_{k}$ be a maximal multi-samplable overlap at inverse-temperature
$\beta_{k}$. Recall that for any $s\in\S$, $q_{k}(s)>0$. By Lemma
\ref{lem:derivative_qb}, $(\beta_{k},q_{k})$ is a solution of (\ref{eq:20051}).

If we define $x_{k}:=\max_{s\in\S}q_{k}(s)$, then the left-hand side
of (\ref{eq:20051}) is bounded from above by
\[
\beta_{k}\Es p(s)x_{k}^{\frac{|p|}{2}-1}-\lambda(s)+\beta_{k}^{2}p(s)x_{k}^{|p|-2}.
\]
Hence, for appropriate $\epsilon=\epsilon(\beta_{c},\Es,p,\lambda)>0$,
we have that $x_{k}>\epsilon$ for all $k$. From (\ref{eq:TAPq>0}),
it is also easy to see that $x_{k}<1-\epsilon$ for all $k$, if $\epsilon$
is small enough.

Suppose that $q_{c}$ is some subsequential limit of $q_{k}$. Note
that (\ref{eq:TAPrep-1}) holds as equality with $(\beta,q)=(\beta_{k},q_{k})$,
since this characterizes multi-samplable overlaps. From the continuity
of free energies in (\ref{eq:TAPrep-1}), it holds as an equality
also for $(\beta,q)=(\beta_{c},q_{c})$. Therefore, $q_{c}$ is a
multi-samplable overlap at $\beta_{c}$ and $q=0$ cannot be maximal.
The lemma thus follows by Lemma \ref{lem:q_nonzero}.\qed

\section{\label{sec:PfThm1}Proof of Theorem \ref{thm:main1}}

Assume that $|p|\geq3$. Suppose that for some $\beta,E>0$ and $q\in(0,1)^{\S}$,
$(\beta,q,E)$ is a solution of the equations in Proposition \ref{prop:eqns}.
From (\ref{eq:II}),
\begin{equation}
\begin{aligned}\beta E\sqrt{\xi(q)}\frac{p(s)}{q(s)} & =-\frac{p(s)}{q(s)}\sum_{t\in\S}\lambda(t)\log(1-q(t)),\\
\beta^{2} & =-\frac{\sum_{t\in\S}\lambda(t)\log(1-q(t))}{\xi(q)\left(1+\sum_{t\in\S}(1-q(t))\frac{p(t)}{q(t)}\right)}.
\end{aligned}
\label{eq:extract}
\end{equation}
 By substituting (\ref{eq:extract}) in (\ref{eq:I}), multiplying
both sides by $1+\sum_{t\in\S}(1-q(t))\frac{p(t)}{q(t)}$ and grouping
together the two terms involving the log we obtain that
\[
\frac{\lambda(s)}{1-q(s)}\left(1+\sum_{t\in\S}(1-q(t))\frac{p(t)}{q(t)}\right)=-\frac{p(s)}{q(s)^{2}}\sum_{t\in\S}\lambda(t)\log(1-q(t)).
\]
This is equivalent to (\ref{eq:qc_formula}). Hence, by Proposition
\ref{prop:eqns}, any maximal multi-samplable overlap at $\beta_{c}$
belongs to $(0,1)^{\S}$ and satisfies (\ref{eq:qc_formula}). The
proof of Part \ref{enu:thm1pt1} of the theorem is completed by the
following lemma which we will prove below.
\begin{lem}
\label{lem:uniq_qc}The system of equations (\ref{eq:qc_formula})
for all $s\in\S$ has at most one solution $q$ in $(0,1)^{\S}$.
\end{lem}

The left-hand side of (\ref{eq:I}) is equal to the left-hand side
of (\ref{eq:II}) times $\frac{p(s)}{q(s)}$. By simple algebra, one
has that
\begin{equation}
\beta^{2}=\frac{1}{\xi(q)}\frac{\lambda(s)}{p(s)}\frac{q(s)^{2}}{1-q(s)}.\label{eq:beta2}
\end{equation}
Part \ref{enu:thm1pt2} of the theorem follows by combining this equality
for $(\beta,q)=(\beta_{c},q_{c})$ with Part \ref{enu:thm1pt1}. Part
\ref{enu:thm1pt3} follows from Part \ref{enu:thm1pt2} and (\ref{eq:II}).
\qed

\subsection{\label{subsec:Proof-of-Lemma_uniqqc}Proof of Lemma \ref{lem:uniq_qc}}

For the convenience of the reader, we recall equations (\ref{eq:qc_formula}),
\begin{equation}
\frac{\lambda(s)}{p(s)}\frac{q(s)^{2}}{1-q(s)}=\Phi(q):=-\frac{\sum_{t\in\S}\lambda(t)\log(1-q(t))}{1+\sum_{t\in\S}p(t)\frac{1-q(t)}{q(t)}}.\label{eq:eq_formula_copy}
\end{equation}

Given $z\geq0$, let $q_{z}(s)$ be the unique solution $x\in[0,1)$
of the equation
\[
f_{s}(x):=\frac{\lambda(s)}{p(s)}\frac{x^{2}}{1-x}=z.
\]
Precisely, we define $q_{z}(s)=f_{s}^{-1}(z)$, where $f_{s}:[0,1)\to[0,\infty)$
is obviously a bijection since it is a strictly increasing function.

Observe that $q\in(0,1)^{\S}$ is a solution of (\ref{eq:eq_formula_copy})
if and only if there exists (a unique) $z>0$ such that $q=q_{z}$
and 
\begin{equation}
z=\Phi(q_{z}).\label{eq:zfixedpt}
\end{equation}
To verify this, note that if $q\in(0,1)^{\S}$ is a solution of (\ref{eq:eq_formula_copy}),
then $f_{s}(q(s))=\Phi(q)$ and with $z=\Phi(q)$,
\[
q_{z}(s)=f_{s}^{-1}(\Phi(q))=q(s),
\]
i.e., $q=q_{z}$. Conversely, if $z$ satisfies (\ref{eq:zfixedpt}),
then 
\[
f_{s}(q_{z}(s))=z=\Phi(q_{z}),
\]
namely, $q_{z}$ solves (\ref{eq:eq_formula_copy}). Hence, to complete
the proof it will be enough to show that the equation (\ref{eq:zfixedpt})
has a unique solution $z>0$.

Equation (\ref{eq:zfixedpt}) is equivalent to 
\[
z\Big(1+\sum_{t\in\S}p(t)\frac{1-q_{z}(t)}{q_{z}(t)}\Big)=-\sum_{t\in\S}\lambda(t)\log(1-q_{z}(t)).
\]
By definition, $z=f_{s}(q_{z}(s))$ and therefore
\[
zp(s)\frac{1-q_{z}(s)}{q_{z}(s)}=\lambda(s)q_{z}(s).
\]
Hence, (\ref{eq:zfixedpt}) is equivalent to $g(z)=0$ for 
\[
g(z):=z+\sum_{t\in\S}\lambda(t)q_{z}(t)+\sum_{t\in\S}\lambda(t)\log(1-q_{z}(t)).
\]
The proof will be completed if we show that the equation $g(z)=0$ has
at most one solution $z>0$. Since $g(0)=0$, it will be enough to
show that $g(z)$ is strictly convex on $[0,\infty)$. 

We wish to compute the derivative $\frac{d}{dz}q_{z}(s)$. To do so,
we first compute the derivative of the corresponding inverse function.
For any $s\in\S$ , 

\[
\frac{d}{dx}f_{s}(x)=\frac{\lambda(s)}{p(s)}\frac{2x-x^{2}}{(1-x)^{2}}.
\]
Therefore for any $z>0$,
\[
\frac{d}{dz}q_{z}(s)=\frac{d}{dz}f_{s}^{-1}(z)=\frac{p(s)}{\lambda(s)}\frac{(1-f_{s}^{-1}(z))^{2}}{2f_{s}^{-1}(z)-f_{s}^{-1}(z)^{2}}=\frac{p(s)}{\lambda(s)}\frac{(1-q_{z}(s))^{2}}{2q_{z}(s)-q_{z}(s)^{2}}.
\]

From this we have that
\begin{align*}
g'(z) & =1+\sum_{t\in\S}p(t)\frac{(1-q_{z}(t))^{2}}{2q_{z}(t)-q_{z}(t)^{2}}-\sum_{t\in\S}p(t)\frac{1-q_{z}(t)}{2q_{z}(t)-q_{z}(t)^{2}}\\
 & =1-\sum_{t\in\S}p(t)\frac{1-q_{z}(t)}{2-q_{z}(t)}=1-\sum_{t\in\S}p(t)\left(1-\frac{1}{2-q_{z}(t)}\right).
\end{align*}
We note that $q_{z}(s)$ is a strictly increasing function of $z>0$.
Hence, $g'(z)$ is also increasing in $z>0$ and therefore $g(z)$
is strictly convex on $[0,\infty)$.\qed

\section{\label{sec:pf_Thm2}Proof of Theorem \ref{thm:main2}}

Let $\beta\geq\beta_{c}$ and suppose $q$ is a maximal multi-samplable
overlap at $\beta$. By Proposition \ref{prop:eqns}, $q\in(0,1)^{\S}$
and for any $s\in\S$,
\[
\frac{\lambda(s)}{p(s)}\frac{1}{\beta\sqrt{\xi(q)}}\frac{q(s)}{1-q(s)}-\beta\sqrt{\xi(q)}\frac{1-q(s)}{q(s)}+\beta\sqrt{\xi(q)}\sum_{t\in\S}p(t)\frac{1-q(t)}{q(t)}=\Es.
\]

Defining
\begin{equation}
\Gamma(s)=\Gamma(s,\beta,q):=\beta\sqrt{\xi(q)}\frac{1-q(s)}{q(s)}\label{eq:xsdef}
\end{equation}
and 
\begin{equation}
\As=\As(\beta,q):=\sum_{t\in\S}p(t)\Gamma(t)=\beta\sqrt{\xi(q)}\sum_{t\in\S}p(t)\frac{1-q(t)}{q(t)},\label{eq:As}
\end{equation}
we may write this as
\begin{equation}
\frac{\lambda(s)}{p(s)}\frac{1}{\Gamma(s)}-\Gamma(s)=\Es-\As.\label{eq:EsAs}
\end{equation}
Note that $\As=\As(\beta,q)$ does not depend on $s\in\S$. We will
also prove below that it is identical for all $\beta$ and $q$ as
above.

By multiplying by $\Gamma(s)$, we obtain from (\ref{eq:EsAs}) a
quadratic equation in $\Gamma(s)$. Solving for $\Gamma(s)$ we have
that it is equal to one of the two values
\[
\frac{-(\Es-\As)\pm\sqrt{(\Es-\As)^{2}+4\frac{\lambda(s)}{p(s)}}}{2}.
\]
Since $q(s)\in(0,1)$, $\Gamma(s)>0$. Hence, we must have that 
\begin{equation}
\Gamma(s)=\frac{-(\Es-\As)+\sqrt{(\Es-\As)^{2}+4\frac{\lambda(s)}{p(s)}}}{2}.\label{eq:xs}
\end{equation}

Plugging this into (\ref{eq:As}) we obtain that
\begin{equation}
\As=\sum_{t\in\S}\frac{p(t)}{2}\left(-(\Es-\As)+\sqrt{(\Es-\As)^{2}+4\frac{\lambda(t)}{p(t)}}\right).\label{eq:As_fixedpt}
\end{equation}
Define, for $y\in\R$,
\begin{equation}
\Theta(y):=y+\sum_{t\in\S}\frac{p(t)}{2}\left(-y+\sqrt{y^{2}+4\frac{\lambda(t)}{p(t)}}\right),\label{eq:Theta}
\end{equation}
so that (\ref{eq:As_fixedpt}) is equivalent to
\begin{equation}
\Theta(\Es-\As)=\Es.\label{eq:AstarEq}
\end{equation}

We will prove the following in Section \ref{subsec:pf_AsUniq}.
\begin{lem}
\label{lem:AsUniq}For any $\beta\geq\beta_{c}$ and any maximal multi-samplable
overlap $q$, $\As=\As(\beta,q)$ satisfies $\Es\geq\As$ and 
\[
\Theta'(\Es-\As)\geq0.
\]
Moreover, $\Theta(y)$ is strictly convex on $\R$, and the value
of $\As$ is therefore identical for all such $\beta$ and $q$. 
\end{lem}

From (\ref{eq:I}), (\ref{eq:bc_formula}) and the definition (\ref{eq:As})
of $\As$, at the critical temperature and with the unique maximal
multi-samplable overlap $q_{c}$,

\[
\As(\beta_{c},q_{c})=\Es-\beta_{c}\sqrt{\xi(q_{c})}=\Es-\sqrt{\Phi(q_{c})}.
\]
Hence, by Lemma \ref{lem:AsUniq}, for any $\beta$ and $q$ as in
the lemma,
\begin{equation}
\As(\beta,q)=\Es-\sqrt{\Phi(q_{c})}.\label{eq:AsVal}
\end{equation}

The value of $\As$ determines $\Gamma(s)$ for all $s\in\S$, but
not $q(s)$ by itself. The pair $(\As,\beta\sqrt{\xi(q)})$, however,
does determine $q(s)$ for all $s\in\S$, since $y\mapsto(1-y)/y$
is bijective on $(0,1)$. Note that
\[
\frac{\Gamma(s)}{\beta\sqrt{\xi(q)}}+1=\frac{1}{q(s)},
\]
and therefore, $y=\beta\sqrt{\xi(q)}$ is a solution of the equation
\begin{equation}
\beta^{2}=\Upsilon(y),\label{eq:UpsilonEq}
\end{equation}
where for $y>0$ we define
\begin{equation}
\begin{aligned}\Upsilon(y) & =y^{2}\prod_{s\in\S}\left(\frac{\Gamma(s)}{y}+1\right)^{p(s)}\\
 & =y^{2}\prod_{s\in\S}\left(\frac{-(\Es-\As)+\sqrt{(\Es-\As)^{2}+4\frac{\lambda(s)}{p(s)}}}{2y}+1\right)^{p(s)}.
\end{aligned}
\label{eq:Upsilon}
\end{equation}

We will prove the following in Section \ref{subsec:pf_ySol}.
\begin{lem}
\label{lem:ySol} Let $\beta>\beta_{c}$ and suppose $q$ is a maximal
multi-samplable overlap at $\beta$. Then there are two positive solutions
to the equation $\Upsilon(y)=\beta^{2}$ and $y=\beta\sqrt{\xi(q)}$
is the larger of the two.
\end{lem}

Combining Lemma \ref{lem:ySol} and (\ref{eq:AsVal}), we have that
if $\beta>\beta_{c}$ and $q$ is a maximal multi-samplable overlap,
then $\beta\sqrt{\xi(q)}$ is equal to the solution $y_{\star}(\beta)$
of (\ref{eq:qb}). The fact that $q$ satisfies (\ref{eq:q}) follows
from (\ref{eq:xs}) and the definition (\ref{eq:xsdef}) of $\Gamma(s)$.
Finally, (\ref{eq:FThm}) follows from the TAP representation (\ref{eq:TAPrep-1})
and (\ref{eq:Onsager}).\qed

\subsection{\label{subsec:pf_AsUniq}Proof of Lemma \ref{lem:AsUniq}}

We have that
\[
\Theta'(y)=1+\sum_{t\in\S}\frac{p(t)}{2}\left(-1+\frac{y}{\sqrt{y^{2}+4\frac{\lambda(t)}{p(t)}}}\right)
\]
and 
\begin{equation}
\begin{aligned}\Theta''(y) & =\sum_{t\in\S}\frac{p(t)}{2}\left(\frac{1}{\sqrt{y^{2}+4\frac{\lambda(t)}{p(t)}}}-\frac{y^{2}}{\left(y^{2}+4\frac{\lambda(t)}{p(t)}\right)^{3/2}}\right)\\
 & =\sum_{t\in\S}\frac{2\lambda(t)}{\left(y^{2}+4\frac{\lambda(t)}{p(t)}\right)^{3/2}}>0.
\end{aligned}
\label{eq:ThetaCvx-1}
\end{equation}
Thus, $\Theta(y)$ is a strictly convex function on $\R$. 

Let $\beta\geq\beta_{c}$ and suppose that $q$ is a maximal multi-samplable
overlap $q$. Recall that $\As=\As(\beta,q)$ satisfies $\Theta(\Es-\As)=\Es$.
From convexity, there are two solutions at most to the equation $\Theta(y)=\Es$.
Note that for $y<0$, $\Theta'(y)<0$. Hence, to complete the proof
of the lemma, it remains to show that 
\begin{equation}
\Theta'(\Es-\As)=1-\sum_{t\in\S}\frac{p(t)}{2}\frac{-(\Es-\As)+\sqrt{(\Es-\As)^{2}+4\frac{\lambda(t)}{p(t)}}}{\sqrt{(\Es-\As)^{2}+4\frac{\lambda(t)}{p(t)}}}\geq0.\label{eq:Theta'star}
\end{equation}

WLOG assume that $\S=\{1,\ldots,|\S|\}$. 
Recall that since $q$ is a maximal
multi-samplable overlap it satisfies \eqref{eq:RS_Gibbs} and, as stated in (\ref{eq:Onsager}),
\[
F(\beta,q)=\frac{1}{2}\beta^{2}\xi_{q}(1).
\]
Namely, $\beta$ is less than or equal to the critical inverse-temperature
of $\xi_{q}(x)$.   
By Lemma 12 of \cite{2ndmomentmulti}, (if $p(s)\geq 2$ for all $s$, Proposition 2 of \cite{2ndmomentmulti} is also sufficient) if we define the function
\[
f_{\beta}(x)=\frac{1}{2}\sum_{s\in\S}\lambda(s)\log(1-x(s)^{2})+\beta^{2}\xi_{q}(x),
\]
then the matrix
\[
\Big(\frac{d}{dx(s)}\frac{d}{dx(t)}f_{\beta}(0)\Big)_{s,t=1}^{|\S|}
\]
is negative semi-definite.

By Proposition \ref{prop:eqns}, $q(s)>0$ for all $s\in\S$. Recall
(\ref{eq:Deltaqp}) and that $\Delta_{p}^{2}=1$ and $\Delta_{p'}^{2}=0$
for any $p'\neq p$. Of course, 
\[
\frac{d}{dx(s)}\frac{d}{dx(t)}\xi_{q}(0)=\xi(q)p(s)(p(t)-\delta_{st})\frac{1-q(s)}{q(s)}\frac{1-q(t)}{q(t)},
\]
where we define
\[
\delta_{st}=\begin{cases}
1 & s=t,\\
0 & s\neq t.
\end{cases}
\]
Hence,
\begin{align*}
\Big(\frac{d}{dx(s)}\frac{d}{dx(t)}f_{\beta}(0)\Big)_{s,t=1}^{|\S|}= & -\bigg(\delta_{st}\Big(\lambda(s)+\beta^{2}\xi(q)p(s)\Big(\frac{1-q(s)}{q(s)}\Big)^{2}\Big)\bigg)_{s,t=1}^{|\S|}\\
 & +\bigg(\beta^{2}\xi(q)p(s)p(t)\frac{1-q(s)}{q(s)}\frac{1-q(t)}{q(t)}\bigg)_{s,t=1}^{|\S|}\\
= & \negthinspace:-M_{1}+M_{2}.
\end{align*}

Defining the column vector
\[
V=\bigg(\beta\sqrt{\xi(q)}p(s)\frac{1-q(s)}{q(s)}\bigg)_{s=1}^{|\S|},
\]
we have that
\[
M_{2}=VV^{T}.
\]

Note that
\[
0\leq\det(M_{1}-VV^{T})=\det(M_{1})\det(I-M_{1}^{-1}VV^{T})
\]
and since $\det(M_{1})>0$,
\[
0\leq\det(I-M_{1}^{-1}VV^{T})=1-V^{T}M_{1}^{-1}V,
\]
where the equality follows from Sylvester's determinant theorem.

Hence,
\begin{equation}
\begin{aligned}1\geq V^{T}M_{1}^{-1}V & =\sum_{s\in\S}\frac{\beta^{2}\xi(q)p(s)^{2}\Big(\frac{1-q(s)}{q(s)}\Big)^{2}}{\lambda(s)+\beta^{2}\xi(q)p(s)\Big(\frac{1-q(s)}{q(s)}\Big)^{2}}\\
 & =\sum_{s\in\S}\frac{p(s)^{2}\Gamma(s)^{2}}{\lambda(s)+p(s)\Gamma(s)^{2}},
\end{aligned}
\label{eq:b}
\end{equation}
where for the second equality we used (\ref{eq:xsdef}).

From (\ref{eq:xs}) and (\ref{eq:EsAs}),
\begin{align*}
 & \frac{p(s)}{2}\frac{-(\Es-\As)+\sqrt{(\Es-\As)^{2}+4\frac{\lambda(s)}{p(s)}}}{\sqrt{(\Es-\As)^{2}+4\frac{\lambda(s)}{p(s)}}}\\
 & =\frac{p(s)\Gamma(s)}{2\Gamma(s)+(\Es-\As)}=\frac{p(s)\Gamma(s)}{\frac{\lambda(s)}{p(s)}\frac{1}{\Gamma(s)}+\Gamma(s)}=\frac{p(s)^{2}\Gamma(s)^{2}}{\lambda(s)+p(s)\Gamma(s)^{2}}.
\end{align*}
Combining this with (\ref{eq:b}) proves (\ref{eq:Theta'star}) and
completes the proof.\qed

\subsection{\label{subsec:pf_ySol}Proof of Lemma \ref{lem:ySol}}

Denote by $q_{c}$ the unique maximal multi-samplable overlap at the
critical inverse-temperature $\beta_{c}$. By Lemma \ref{lem:derivative_bc},
\[
\beta_{c}\xi(1)=\Es\sqrt{\xi(q_{c})}+\beta_{c}\left(\xi(1)-\xi(q_{c})-\xi(q_{c})\sum_{t\in\S}p(t)\frac{1-q_{c}(t)}{q_{c}(t)}\right).
\]
Using the definition of $\As$ in (\ref{eq:As}), we have that 
\[
\beta_{c}\sqrt{\xi(q_{c})}=\Es-\As.
\]

Now, let $\beta>\beta_{c}$ and suppose that $q$ is a maximal multi-samplable
overlap. With the same $q$, for any $\beta'$, recall that the TAP
representation (\ref{eq:TAPrep-1}) holds as an inequality,
\[
F(\beta')\geq\beta'\sqrt{\xi(q)}\Es+\frac{1}{2}\sum_{s\in\S}\lambda(s)\log(1-q(s))+F(\beta',q).
\]
For $\beta'=\beta$, the above becomes an equality. Since $\beta\mapsto F(\beta)$
is convex, the derivative from the right $\frac{d}{d\beta}^{+}F(\beta)$
exists. Hence,
\[
\frac{d}{d\beta}^{+}F(\beta)\geq\sqrt{\xi(q)}\Es+\frac{d}{d\beta}F(\beta,q).
\]
As before, using Lemma \ref{lem:derivative_bc} and (\ref{eq:As}),
we obtain that
\begin{equation}
\begin{aligned}\beta\sqrt{\xi(q)} & \geq\Es-\As+\frac{1}{\sqrt{\xi(q)}}\Big(\beta\xi(1)-\frac{d}{d\beta}^{+}F(\beta)\Big)\\
 & =\beta_{c}\sqrt{\xi(q_{c})}+\frac{1}{\sqrt{\xi(q)}}\Big(\beta\xi(1)-\frac{d}{d\beta}^{+}F(\beta)\Big),
\end{aligned}
\label{eq:bsqrtxi}
\end{equation}
where we recall that by Proposition \ref{prop:eqns}, $q(s)>0$ and
therefore $\xi(q)>0$ and that by Lemma \ref{subsec:pf_AsUniq}, $\As$
does not depend on $\beta$ or $q$.

Note that, in distribution, for any $\beta<\beta'$,
\[
\beta'H_{N}(\bs)=\beta H_{N}^{1}(\bs)+\sqrt{\beta'^{2}-\beta^{2}}H_{N}^{2}(\bs),
\]
where $H_{N}^{i}(\bs)$ are i.i.d. copies of the Hamiltonian $H_{N}(\bs)$.
By applying Jensen's inequality conditional on $H_{N}^{1}(\bs)$,
one sees that
\begin{align*}
\frac{1}{N}\E\log\int_{S_{N}}e^{\beta'H_{N}(\bs)}d\mu(\bs) & =\frac{1}{N}\E\log\int_{S_{N}}e^{\beta H_{N}^{1}(\bs)+\sqrt{\beta'^{2}-\beta^{2}}H_{N}^{2}(\bs)}d\mu(\bs)\\
 & \leq\frac{1}{N}(\beta'^{2}-\beta^{2})\E H_{N}^{2}(\bs)+\frac{1}{N}\E\log\int_{S_{N}}e^{\beta H_{N}^{1}(\bs)}d\mu(\bs)\\
 & =\frac{1}{2}(\beta'^{2}-\beta^{2})\xi(1)+\frac{1}{N}\E\log\int_{S_{N}}e^{\beta H_{N}(\bs)}d\mu(\bs).
\end{align*}
Hence, 
\[
\frac{d}{d\beta}^{+}F(\beta)\leq\beta\xi(1).
\]
Combined with (\ref{eq:bsqrtxi}), this gives that for any $\beta\geq\beta_{c}$,
\begin{equation}
\beta\sqrt{\xi(q)}\geq\beta_{c}\sqrt{\xi(q_{c})}.\label{eq:bxibd}
\end{equation}

Recall the definition (\ref{eq:Upsilon}) of $\Upsilon(y)$. For $y>0$,
\begin{equation}
\Upsilon'(y)=\Upsilon(y)\left(\frac{2}{y}-\sum_{s\in\S}\frac{p(s)}{\frac{\Gamma(s)}{y}+1}\frac{\Gamma(s)}{y^{2}}\right)=\frac{\Upsilon(y)}{y}\left(2-\sum_{s\in\S}\frac{p(s)\Gamma(s)}{\Gamma(s)+y}\right).\label{eq:Upsilon'}
\end{equation}
(Note that since $\As$ does not depend on $\beta$ or $q$, so does
$\Gamma(s)$.) Hence, for some $y_{0}>0$, $\Upsilon(y)$ strictly
decreases on $(0,y_{0}]$ and strictly increases on $[y_{0},\infty)$.
Clearly, 
\[
\lim_{y\searrow0}\Upsilon(y)=\lim_{y\to\infty}\Upsilon(y)=\infty.
\]

Recall that $y=\beta_{c}\sqrt{\xi(q_{c})}$ is a solution to the equation
$\Upsilon(y)=\beta_{c}^{2}$, (see (\ref{eq:UpsilonEq})) and thus
$\Upsilon(y_{0})\leq\beta_{c}^{2}$. If $\Upsilon(y_{0})<\beta_{c}^{2}$,
then there are exactly two positive solutions to the same equation
$\Upsilon(y)=\beta_{c}^{2}$. Denote those solutions by $y_{c}^{(1)}$
and $y_{c}^{(2)}$ and assume that $y_{c}^{(1)}\leq y_{c}^{(2)}$.
If $\Upsilon(y_{0})=\beta_{c}^{2}$ then there is one solution and
it is equal to $y_{0}$. In this case we will assume that $y_{c}^{(1)}=y_{c}^{(2)}=y_{0}$.

For $\beta>\beta_{c}$, the equation $\Upsilon(y)=\beta^{2}$ has
exactly two positive solutions, and by (\ref{eq:UpsilonEq}), $y=\beta\sqrt{\xi(q)}$
is one of those solutions (where $q$ is some maximal multi-samplable
overlap at $\beta$). From the properties of $\Upsilon(y)$ we mentioned
above, one of those two solutions is strictly less than $y_{c}^{(1)}$
and the other is strictly greater than $y_{c}^{(2)}$. By (\ref{eq:bxibd})
and since $\beta_{c}\sqrt{\xi(q_{c})}\in\{y_{c}^{(1)},\,y_{c}^{(2)}\}$,
we therefore have that $\beta\sqrt{\xi(q)}$ is the larger of the
two solutions to the equation $\Upsilon(y)=\beta^{2}$.\qed

\section*{Appendix: the case $|p|=2$}

In this appendix we analyze the free energy of the pure $p$-spin
multi-species (bipartite) model with $|p|=2$, assuming the convergence
condition from (\ref{eq:convergence}).

Suppose that $\S=\{s,t\}$ and $p(s)=p(t)=1$. The corresponding Hamiltonian
can be written as
\[
H_{N}(\bs)=\sqrt{\frac{N}{N_{s}N_{t}}}\sum_{i\in I_{s},\,j\in I_{t}}J_{ij}\sigma_{i}\sigma_{j}.
\]
In distribution, 
\[
\frac{1}{N}\max_{\bs\in S_{N}}H_{N}(\bs)=\frac{1}{\sqrt{N}}\max_{\|\bx\|=\|\by\|=1}\bx^{T}M\by=\sqrt{\frac{1}{N}\nu_{\max}(M^{T}M)},
\]
where $M$ is a $|I_{s}|\times|I_{t}|$ matrix whose elements are
i.i.d standard normal variables, $\nu_{\max}(A)$ is the maximal eigenvalue
of $A$, and $\bx$ and $\by$ are real column vectors of length $|I_{s}|$
and $|I_{t}|$, respectively.

Note that $M^{T}M$ is a Wishart matrix. The maximal eigenvalue of
a Wishart matrix was studied in \cite{GemanWishart} where it was
shown that, when properly normalized, it converges a.s. to an explicit
deterministic limit. Combining the latter with the Borell-TIS inequality,
we have that 
\begin{equation}
\Es=\sqrt{\lambda(s)}+\sqrt{\lambda(t)}.\label{eq:Es_p=00003D2}
\end{equation}

The fact that $|p|\geq3$ was used in the proof of Proposition \ref{prop:eqns}
in Lemma \ref{lem:q>0}, which  states that at $\beta_{c}$ all the
elements $q(s)$ of a maximal multi-samplable are strictly positive
(which we will soon see is not the case for $|p|=2$). However, all
arguments in same the proof concerning $\beta>\beta_{c}$ apply as
they are also for $|p|=2$. In particular, for any $\beta>\beta_{c}$,
if $q$ is a maximal multi-samplable overlap then, for all $s\in\S$,
$q(s)>0$ and $(\beta,q,\Es)$ solves (\ref{eq:I}). Further, with
$\As=\As(\beta,q)$ as defined in (\ref{eq:As}), by the same argument
as in Section \ref{sec:pf_Thm2} we have that $\As$ satisfies (\ref{eq:As_fixedpt}).
In the current setting with $p(s)=p(t)=1$, by rearranging the latter
equation we obtain that
\[
\Es=\frac{1}{2}\sqrt{(\Es-\As)^{2}+4\lambda(s)}+\frac{1}{2}\sqrt{(\Es-\As)^{2}+4\lambda(t)},
\]
which in light of (\ref{eq:Es_p=00003D2}) implies that, for all $\beta>\beta_{c}$
and $q$ as above,
\[
\Es=\As.
\]

Let $\beta>\beta_{c}$ and suppose that $q$ is a maximal multi-samplable
overlap. Recall that $y=\beta\sqrt{\xi(q)}$ solves (\ref{eq:UpsilonEq}),
namely,
\begin{align*}
\beta^{2} & =y^{2}\left(\frac{\sqrt{\lambda(s)}}{y}+1\right)\left(\frac{\sqrt{\lambda(t)}}{y}+1\right)\\
 & =y^{2}+y(\sqrt{\lambda(s)}+\sqrt{\lambda(t)})+\sqrt{\lambda(s)\lambda(t)}.
\end{align*}
The two solutions for the latter quadratic equation are
\begin{align*}
y_{\pm} & =-\frac{1}{2}\left(\sqrt{\lambda(s)}+\sqrt{\lambda(t)}\right)\pm\frac{1}{2}\sqrt{\left(\sqrt{\lambda(s)}+\sqrt{\lambda(t)}\right)^{2}+4(\beta^{2}-\sqrt{\lambda(s)\lambda(t)})}\\
 & =-\frac{1}{2}\left(\sqrt{\lambda(s)}+\sqrt{\lambda(t)}\right)\pm\frac{1}{2}\sqrt{\left(\sqrt{\lambda(s)}-\sqrt{\lambda(t)}\right)^{2}+4\beta^{2}}.
\end{align*}
Since $\beta\sqrt{\xi(q)}>0,$ we have that $\beta\sqrt{\xi(q)}=y_{+}$.

From (\ref{eq:xs}), $\Gamma(s)=\sqrt{\lambda(s)}$ and $\Gamma(t)=\sqrt{\lambda(t)}$.
From the definition (\ref{eq:xsdef}) of $\Gamma(s)$,
\begin{align*}
1-q(s) & =\frac{\Gamma(s)}{\Gamma(s)+\beta\sqrt{\xi(q)}}\\
 & =\frac{2\sqrt{\lambda(s)}}{\sqrt{\lambda(s)}-\sqrt{\lambda(t)}+\sqrt{\left(\sqrt{\lambda(s)}-\sqrt{\lambda(t)}\right)^{2}+4\beta^{2}}}
\end{align*}
and
\begin{equation}
\begin{aligned}q(s) & =\frac{\beta\sqrt{\xi(q)}}{\Gamma(s)+\beta\sqrt{\xi(q)}}\\
 & =\frac{-\sqrt{\lambda(s)}-\sqrt{\lambda(t)}+\sqrt{\left(\sqrt{\lambda(s)}-\sqrt{\lambda(t)}\right)^{2}+4\beta^{2}}}{\sqrt{\lambda(s)}-\sqrt{\lambda(t)}+\sqrt{\left(\sqrt{\lambda(s)}-\sqrt{\lambda(t)}\right)^{2}+4\beta^{2}}}
\end{aligned}
\label{eq:qp=00003D2}
\end{equation}
and similarly for $t$.

Using the fact that $p(s)=p(t)=1$,
\[
\xi_{q}(1)=\left(1-q(s)\right)\left(1-q(t)\right)=\frac{1}{\beta^{2}}\sqrt{\lambda(s)\lambda(t)}.
\]

Combining the above with the TAP representation (\ref{eq:TAPrep-1})
and the correction (\ref{eq:Onsager}), after some algebra we obtain
that
\[
\begin{aligned}F(\beta) & =\frac{1}{2}\left(-1-\sqrt{\lambda(s)\lambda(t)}+\left(\sqrt{\lambda(s)}+\sqrt{\lambda(t)}\right)\sqrt{\left(\sqrt{\lambda(s)}-\sqrt{\lambda(t)}\right)^{2}+4\beta^{2}}\right)\\
 & -\frac{\lambda(s)-\lambda(t)}{4}\log\left(\frac{\sqrt{\lambda(s)}-\sqrt{\lambda(t)}+\sqrt{\left(\sqrt{\lambda(s)}-\sqrt{\lambda(t)}\right)^{2}+4\beta^{2}}}{\sqrt{\lambda(t)}-\sqrt{\lambda(s)}+\sqrt{\left(\sqrt{\lambda(s)}-\sqrt{\lambda(t)}\right)^{2}+4\beta^{2}}}\right)\\
 & -\frac{1}{2}\log\beta+\frac{1}{4}\lambda(s)\log(\lambda(s))+\frac{1}{4}\lambda(t)\log(\lambda(t)).
\end{aligned}
\]
This gives the free energy for any $\beta>\beta_{c}$. Reassuringly,
this coincides with the formula from Theorem 2.1 of \cite{BaikLeeBipartite},
if we replaces $\beta$ by $\sqrt{\lambda(s)\lambda(t)}\beta$ to
account for the different normalization of the Hamiltonian in the
same paper. 

Recall that by the proof of Lemma \ref{lem:derivative_bc}, $F'(\beta_{c})=\beta_{c}$.
Hence, $\beta=\beta_{c}$ satisfies $F'(\beta)-\beta=0$, where here
$F'(\beta)$ is the derivative of formula for $F(\beta)$ above. One
can check that the equation $F'(\beta)-\beta=0$ is equivalent to
\begin{equation}
\kappa(\beta):=\left(\sqrt{\lambda(s)}+\sqrt{\lambda(t)}\right)\sqrt{\left(\sqrt{\lambda(s)}-\sqrt{\lambda(t)}\right)^{2}+4\beta^{2}}-1-2\beta^{2}=0.\label{eq:kappa}
\end{equation}
Note that $\beta=(\lambda(s)\lambda(t))^{\frac{1}{4}}$ solves this
equation. Also, since
\[
\kappa'(\beta)=4\beta\left(\frac{\left(\sqrt{\lambda(s)}+\sqrt{\lambda(t)}\right)}{\sqrt{\left(\sqrt{\lambda(s)}-\sqrt{\lambda(t)}\right)^{2}+4\beta^{2}}}-1\right),
\]
the point $\beta=(\lambda(s)\lambda(t))^{\frac{1}{4}}$ is in fact
also the global maximizer of $\kappa(\beta)$ on $[0,\infty)$. Hence,
$\beta_{c}=(\lambda(s)\lambda(t))^{\frac{1}{4}}$.

Lastly, note that by (\ref{eq:qp=00003D2}), at $\beta_{c}$ the maximal
multi-samplable overlap is $q\equiv0$. 

\bibliographystyle{plain}
\bibliography{master}

\end{document}